\documentclass[11pt]{elsarticle}

\usepackage{comment}
\usepackage{geometry}
\usepackage{natbib}
\usepackage{url}
\usepackage{soul}

\setcitestyle{round,authoryear}
\usepackage{hyperref}
\usepackage[english]{babel}
\addto\extrasenglish{
    
}
\usepackage{array}
\usepackage{multirow}
\usepackage{tabularx}
\usepackage{bm}
\usepackage{amssymb}
\usepackage{booktabs}
\usepackage{caption}
\usepackage{longtable}
\usepackage{amsfonts,amsmath}
\usepackage{xcolor}
\usepackage[linesnumbered]{algorithm2e}
\allowdisplaybreaks
\usepackage{bbm}
\usepackage{tikz}
\usepackage{float}
\usepackage{subfig}
\usepackage{amsmath}
\usepackage[normalem]{ulem}
\usetikzlibrary{shapes.misc, positioning}
\usetikzlibrary{arrows.meta}
\usepackage{mdframed}
\usepackage{lscape}
\usepackage{listings}

\newcounter{cdcounter}  
\definecolor{lightgray}{gray}{0.9}

\newcommand{\btHL}[1]{\colorbox{lightgray}{#1}}

\makeatletter
\def\ps@pprintTitle{%
  \let\@oddhead\@empty
  \let\@evenhead\@empty
  \def\@oddfoot{\reset@font\hfil\thepage\hfil}
  \let\@evenfoot\@oddfoot
}
\makeatother

\begin{document}

\begin{frontmatter}



\title{An LLM-powered MILP modelling engine for workforce scheduling guided by expert knowledge}




\author[label1,label2]{Qingyang Li\corref{cor1}}
\ead{ql5@student.unimelb.edu.au}
\author[label1,label2]{Lele Zhang}
\ead{lele.zhang@unimelb.edu.au}
\author[label3]{Vicky Mak-Hau}
\ead{vicky.mak@deakin.edu.au}

\affiliation[label1]{organization={The University of Melbourne},
              city={Melbourne},
              postcode={3010},
             state={Victoria},
             country={Australia}}
\affiliation[label2]{organization={ARC Training Centre in Optimisation Technologies, Integrated Methodologies, and Applications (OPTIMA)},
             country={Australia}}
\affiliation[label3]{organization={Deakin University},
              city={Burwood},
              postcode={3125},
             state={Victoria},
             country={Australia}}
\cortext[cor1]{Corresponding author.}

\begin{abstract}

Formulating mathematical models from real-world decision problems is a core task in Operational Research, yet it typically requires considerable human expertise and effort, limiting practical application.
Recent advances in large language models (LLMs) have sparked interest in automating this process from natural language descriptions. 
However, challenges including limited modelling expertise, dependence on large-scale training data, and hallucination affect the reliable application of LLMs in optimisation modelling.
To address these challenges, we propose SMILO, an expert-knowledge-driven framework that integrates optimisation modelling expertise with LLMs to generate mixed-integer linear programming models. 
SMILO uses a three-stage architecture built on reusable modelling graphs and associated resources: identifying relevant modelling components, extracting instance-specific information using LLMs, and constructing models through expert-defined templates.
This modular architecture separates information extraction from formula generation, enhancing modelling accuracy, transparency, and reproducibility.
We demonstrate the implementation of our problem-type-specific modelling framework using workforce scheduling problems spanning manufacturing, logistic, and service operations as illustrative cases.
Experiments show that SMILO consistently generates correct models in 90\% of test instances across five trials, outperforming one-step LLM baselines by at least 35\%. 
This work offers a generalisable paradigm for integrating LLMs with expert knowledge across diverse decision-making contexts, advancing automation in optimisation modelling. 

\end{abstract}



\end{frontmatter}

\section{Introduction}
 
Mathematical optimisation modelling plays a fundamental role in effective decision-making across a wide range of domains, including 
logistics, inventory management, production planning, and workforce scheduling \citep{chen_dynamic_2022, liu_food_2024, chen_constructing_2024}. 
A central challenge in applying optimisation lies in the accurate formulation of real-world problems as mathematical models, typically in the form of mixed-integer linear programming (MILP) models.
While advancements in optimisation solvers have improved computation efficiency, the task of translating problem descriptions in natural language into MILP models remains non-trivial.
The modelling process involves the definition of decision variables, identification of constraints, and construction of an appropriate objective function. The process is inherently time-consuming and error-prone, requiring substantial domain knowledge and expertise to accurately convert complex operational requirements into mathematical terms.
This modelling step constitutes a critical bottleneck, limiting broader adoption and application of optimisation techniques, especially for organisations lacking access to dedicated optimisation expertise.

Recent breakthroughs in large language models (LLMs) have inspired substantial research efforts exploring their potential to automate the process of optimisation model formulation from natural language descriptions. LLMs have showcased remarkable proficiency in tasks such as code generation, text comprehension, and even preliminary attempts at mathematical model formulation \citep{jackson_generative_2024, boudribila_intelligent_2025}.  
Nonetheless, their effectiveness remains limited when applied to even moderately complex MILP modelling tasks.
Recent studies have highlighted several fundamental challenges associated with the use of LLMs for MILP model generation \citep{li_synthesizing_2023, ahmaditeshnizi_optimus-03_2024,huang_orlm_2025}:
\begin{itemize}
    \item Challenge 1: \textbf{Limited modelling expertise.} Despite being trained on vast amounts of general textual data, LLMs lack specialised knowledge and expertise in mathematical optimisation. As a result, they often produce incomplete or incorrect models for complex real-world decision-making problems, misdefining decision variables, omitting implicit constraints, miscomputing parameters, or applying inappropriate modelling techniques.

    \item Challenge 2: \textbf{Insufficient large-scale, high-quality training data.} The scarcity of high-quality training data comprising natural language problem descriptions and their corresponding optimisation models constrains the effectiveness of data-driven methods.
    
    \item Challenge 3: \textbf{Limited complexity in benchmark datasets.} While existing benchmark datasets span a wide range of domains such as retail, manufacturing, and agriculture, most instances are simple linear programming problems with structurally similar formulations. Classic combinatorial optimisation problems such as shift scheduling and vehicle routing are included but with few instances per problem type, offering limited variation in modelling complexity. 
    
    \item Challenge 4: \textbf{Hallucination.} LLMs are prone to generating outputs, such as constraints in MILP models, that appear plausible but are in fact incorrect or fabricated, and this affects the reliability of the models generated by the LLMs.
\end{itemize}

To address these challenges, we propose an expert-knowledge-driven MILP modelling framework powered by LLMs, and we name it as SMILO. 
This customised modelling engine integrates structured expert knowledge into a modular and interpretable modelling pipeline.
The framework is underpinned by \textit{problem-type-specific modelling graphs} (MGs), which explicitly organise core modelling components, such as decision variables, constraints, and objective functions, and capture their logical interdependencies. 

The modelling process follows a three-stage pipeline. 
First, sentences from the problem description are matched with nodes in the MG to identify relevant modelling components.
Second, an LLM guided by carefully tailored task-specific templates extracts instance-specific numerical values for the identified components. These templates are designed to guide the LLM's attention to precise information needs, thereby reducing hallucination and improving output consistency.
Third, MILP model construction is carried out using predefined templates that automatically incorporate the extracted values. Parameter calculations and expression generation are programmatically performed using expert-defined rules to ensure mathematical correctness and consistency.

Unlike data-driven approaches that typically aim to enhance the modelling capabilities of LLMs through training and fine-tuning, our method does not delegate core modelling responsibilities to the LLM. Instead, it decomposes the modelling process into structured subtasks, with LLMs assisting in natural language understanding and information extraction, while mathematical expressions are constructed through expert-defined rules and templates. 
This separation of roles and the modular design of the framework enhance modelling accuracy, control, and interpretability, while significantly reducing hallucination-related errors (addressing both Challenge 1 and Challenge 4). 
Furthermore, our framework adopts a knowledge-driven strategy that avoids reliance on large training datasets, relying instead on structured expert knowledge encoded in reusable model components (addressing Challenge 2).
In addition, the modular framework is designed to be \textit{model-neutral}; the LLM component can be seamlessly replaced with any open-source or closed-source models without the need for additional training or significant system adjustments. This flexibility is particularly valuable for domains with strict data privacy requirements or limited computational resources.

To validate the feasibility and effectiveness of our approach, we present a proof-of-concept application targeting two classical yet complex workforce scheduling problems: \textit{shift scheduling} and \textit{days-off scheduling}. These problems exhibit considerable modelling complexity, involving features such as shift coverage, overtime rules, and break allocation. 
In such cases, a direct translation of natural language problem descriptions into mathematical expressions is often infeasible.
While we implement and evaluate our framework on these scheduling problems as illustrative cases, the underlying framework demonstrates a generalisable paradigm broadly applicable to a wide range of MILP problems beyond workforce scheduling. This structured paradigm provides a practical and robust template for effectively leveraging LLM's capability in OR applications.

The primary contributions of our work include:
\begin{itemize}
    \item \textbf{Expert-knowledge-driven modelling engine}: We develop a novel MILP modelling engine that leverages LLMs in a structured, modular, and highly interpretable manner, enabling accurate, reliable, and reproducible model generation from natural language problem descriptions. This modelling engine follows a three-stage framework: component identification, information extraction, and model generation. The framework delegates language understanding tasks to LLMs while retaining precise control of the mathematical formulation through expert knowledge.

    \item \textbf{Problem-type-specific modelling graphs (MGs)}: We introduce and implement MGs as structured representations of reusable modelling components and their interdependencies. These graphs serve as the backbone for guiding component identification, information extraction, and formula assembly. The MGs are modular, extensible, and reusable across similar problem types. 
    We construct a set of reusable modelling resources aligned with the MGs, including predefined formula templates, prompt templates, modelling component descriptions, parameter computation rules, and symbol definition rules, which support systematic and interpretable model generation. These resources enhance transparency, consistency, and traceability throughout the modelling pipeline.

    \item \textbf{Task-based prompting strategy}: We design a task-based prompting strategy that interfaces with the MG. For each activated node in the graph, a corresponding specialised prompt template extracts instance-specific information via the LLM. This approach allows the framework to handle complex and ambiguous problem descriptions and enhances robustness by clearly defining task boundaries and output formats.
    
    \item \textbf{Benchmark dataset for workforce scheduling}: We build a benchmark dataset focused on two representative workforce scheduling problems. The dataset captures varying levels of complexity and structural richness, enabling a more realistic and rigorous evaluation of automated modelling tools beyond toy-sized instances (addressing Challenge 3).
    
    \item \textbf{Comprehensive experimental evaluations}: We evaluate our framework against standard one-step model generation approaches. Across five independent trials, our framework consistently generated correct MILP models for 90\% of instances, while the best-performing baseline achieved no more than 55\% correctness in any single trial. We also provide a detailed error analysis comparing the two approaches, offering valuable insights into common failures and potential areas for future improvement in LLM-assisted optimisation modelling. 
    
    \item \textbf{General and LLM-agnostic paradigm}: We propose a generalisable, LLM-agnostic paradigm for optimisation modelling that requires no fine-tuning or extensive training data. Our framework is readily compatible with various LLMs and extendable to a broad range of optimisation problems, making it an adaptable and scalable solution for real-world decision-support applications.
\end{itemize}

The remainder of this paper is structured as follows. Section \ref{sec:literature} reviews the literature on LLM-based optimisation modelling approaches and workforce scheduling problems. Section \ref{sec:methodology} details the proposed framework. Section \ref{sec:experimental} describes the experimental setup and dataset and presents the experimental results. Finally, Section \ref{sec:conclusion} outlines some conclusions and discusses future work.

\section{Literature Review} \label{sec:literature}

\subsection{LLM-based approaches for automated model generation}

The automation of mathematical programming model generation has gained increasing attention, particularly with the rise of LLMs and their ability to process complex problem descriptions \citep{fan_artificial_2024}. 
While optimisation modelling has traditionally been a manual and expertise-driven process, recent developments have explored various methods to bridge the gap between natural language problem descriptions and formal mathematical representations \citep{ramamonjison22augmenting}. The integration of natural language processing (NLP) techniques, particularly LLMs, into mathematical programming model generation has introduced a new paradigm in automated model generation. 

LLMs, such as GPT, Gemini, PaLM and Llama, have demonstrated strong capabilities in semantic understanding, structured information extraction, and mathematical reasoning, making them suitable for processing optimisation problems described in natural language \citep{fosso_wamba_chatgpt_2024, makatura_how_2023, jackson_natural_2024}. However, despite these advancements, LLMs still exhibit limitations in numerical consistency, logical reasoning, and constraint completeness, which necessitate additional mechanisms to ensure the correctness of the generated models.
Without providing any examples or specific knowledge of MILP, models generated by LLMs frequently contain errors and fail to accurately capture the requirements specified in the problem description, particularly for problems beyond textbook-level difficulty \citep{ramamonjison23nl4opt, li_synthesizing_2023, fan_artificial_2024}.

Existing research has shown that LLMs can be employed to extract variables, constraints, and objective functions from textual inputs, using multi-stage frameworks, fine-tuning, few-shot learning, prompt engineering, and continued pretraining to improve accuracy \citep{li_synthesizing_2023, fan_artificial_2024, li_nl2or_2024, tang_orlm_2024, huang_orlm_2025}.
Other studies improve the correctness and relevance of generated models by incorporating additional validation mechanisms and agent-based adjustment approaches. Multi-agent collaboration enables model generation to be an iterative process and, in some cases, even interactive with users \citep{gaddam_agentmilo_2025, abdullin24synthetic, ahmaditeshnizi23optimus, ahmaditeshnizi_optimus-03_2024, li23large, xiao_chain--experts_2023}.  
Furthermore, a chatbot has been developed to interpret optimisation models to users and diagnose sources of infeasibility \citep{chen_diagnosing_2023}.  
Regarding the final output, existing research focuses on different aspects, including the mathematical model itself and the solution obtained from solving the model \citep{ahmed_lm4opt_2024, ramamonjison23nl4opt, tsouros23holy, ahmaditeshnizi_optimus-03_2024, li23large}.
Mathematical models can be represented in two forms: the canonical form, which expresses the model using matrix notation, and the structured representation, which explicitly defines sets, parameters, variables, constraints, and the objective function. The structured representation offers greater interpretability.

To address LLMs' limitations in domain-specific knowledge, promising approaches include integrating them with task-relevant specialised data, such as Retrieval-Augmented Generation (RAG) \citep{gao_retrieval-augmented_2024}, and incorporating structured knowledge representations \citep{kau_combining_2024, kosasih_towards_2024}.
OptiMUS-0.3 uses RAG to retrieve examples of constraints from their dataset that are similar to the constraint description that is currently being formulated, incorporating them and their formulas into the prompt to enhance LLMs' performance in mathematical modelling. The provided examples can illustrate modelling techniques to LLMs, such as the Big-M technique \citep{ahmaditeshnizi_optimus-03_2024}.

To the best of our knowledge, structured knowledge representations have not yet been integrated into LLM-based automated modelling approaches. In this work, we propose problem-type-specific modelling graphs (MGs) that represent structured expert knowledge of optimisation modelling. These MGs help ensure that essential modelling components are considered and that no critical constraints are overlooked when processing a given problem instance.
Most existing studies focus on LLMs' intrinsic ability to formulate mathematical models. Our research seeks to separate the natural language processing and mathematical formulation stages by using MGs to assist LLMs in identifying key components and extracting relevant information from the problem description. The extracted information is then mapped onto predefined mathematical templates to construct complete MILP models. This approach allows LLMs to focus solely on tasks they are inherently proficient at. This integration ensures both adaptability to complex modelling scenarios and reliability in mathematical model correctness.

\subsection{Workforce scheduling problems}

Workforce scheduling is a fundamental optimisation problem in Operational Research, extensively studied due to its practical significance in various industries such as healthcare, call centres, manufacturing, and banking \citep{van_den_bergh_personnel_2013, mutlu_co-availability_2015}. 
Workforce scheduling problems exhibit considerable diversity depending on the industry, constraints, and objectives, and they are generally categorised into three primary subtypes: shift scheduling, days-off scheduling, and tour scheduling \citep{morris_simple_1983}.
This subsection does not attempt a comprehensive review of workforce scheduling. Instead, it focuses on the two subtypes examined in this study, shift scheduling and days-off scheduling, by reviewing key modelling approaches that were considered during the development of our modelling graphs. Tour scheduling falls outside the scope of this work.

Shift scheduling, also known as time-of-day scheduling, involves determining the start time of shifts, their durations, and how to allocate employees across a planning horizon to meet demand requirements, while incorporating break placements to comply with labour regulations. The shift scheduling problem is particularly crucial in industries where labour demand fluctuates throughout the day.
The shift scheduling problem has been extensively studied in the literature since its introduction by \cite{edie_traffic_1954}. Early Integer Programming (IP) formulations were based on the set-covering formulations proposed by \cite{dantzig_letter_1954}, which define variables by enumerating all possible shifts based on different combinations of start times, durations, and break placements. However, these approaches often result in large-scale IP models with a vast number of variables, making them computationally impractical for real-world applications.  

Subsequent research has sought to refine IP formulations to improve computational efficiency. \cite{gaballa_telephone_1979} developed an implicit IP formulation that introduced additional variables to model lunch break windows. However, their approach required a larger number of integer variables and constraints compared to the equivalent set-covering formulation, making it less efficient. \cite{bechtold_implicit_1990} introduced integer variables to represent the number of employees starting their breaks within specified planning periods, implicitly matching breaks to explicitly represented shifts. A double implicit shift scheduling model was proposed as an extension of \cite{bechtold_implicit_1990}, implicitly matching meal breaks to implicitly represented shifts \citep{thompson_improved_1995}. However, their models were restricted to a single meal break.  

For solving large-scale shift scheduling problems with multiple breaks and multiple break windows, \cite{aykin_optimal_1996} proposed a new implicit IP model by introducing a set of break variables for each shift-break type combination, significantly reducing the number of decision variables and improving computational performance. This approach leads to a substantially smaller number of variables and improved solvability for large instances involving multiple breaks and break windows.  

Days-off scheduling, also known as days-of-week scheduling, focuses on assigning employees' work and rest days over a given scheduling horizon (e.g., a week or a month) to ensure that labour demands are met while complying with contractual agreements and work-life balance considerations. 
A fundamental problem in this category is the $(k, m)$-cyclic workforce scheduling problem, where each worker is assigned a repeating work schedule over an $m$-day cycle, ensuring that they work for $k$ consecutive days and then take $m - k$ days off. It can be formulated by the set-covering model \citep{bartholdi_cyclic_1980}. 

Tour scheduling integrates both shift and days-off scheduling to construct employee work schedules that span multiple days or weeks. This approach is more comprehensive as it considers both intra-day shift assignments and inter-day work patterns, ensuring a globally optimised workforce allocation \citep{alfares_survey_2004}. 
An IP model was developed to address tour scheduling, incorporating additional constraints such as requirements for full-time and part-time workers, consecutive days off, and variable daily shift start times \citep{bard_staff_2003}. 
Researchers have also explored hierarchical scheduling frameworks based on IP models, where high-level workforce allocation decisions are made first, followed by detailed assignments at a finer granularity \citep{turker_integrated_2018}. 


\section{Methodology}
\label{sec:methodology}

\subsection{Framework overview}

The proposed framework automates the transformation of unstructured natural language problem descriptions into mixed-integer linear programming models through a structured three-stage process. First, relevant modelling components are identified by mapping the textual description to nodes in a predefined modelling graph. Second, for the activated nodes, the framework employs task-specific prompts to extract instance-specific numerical values and contextual details using a large language model. Third, the extracted information is processed into sets, parameters, and decision variables, which are then assembled into a complete MILP model (in LaTeX format) using predefined formula templates and rule-based logic.

\begin{figure}[!ht]
  \centering
  \includegraphics[width=0.9\linewidth]{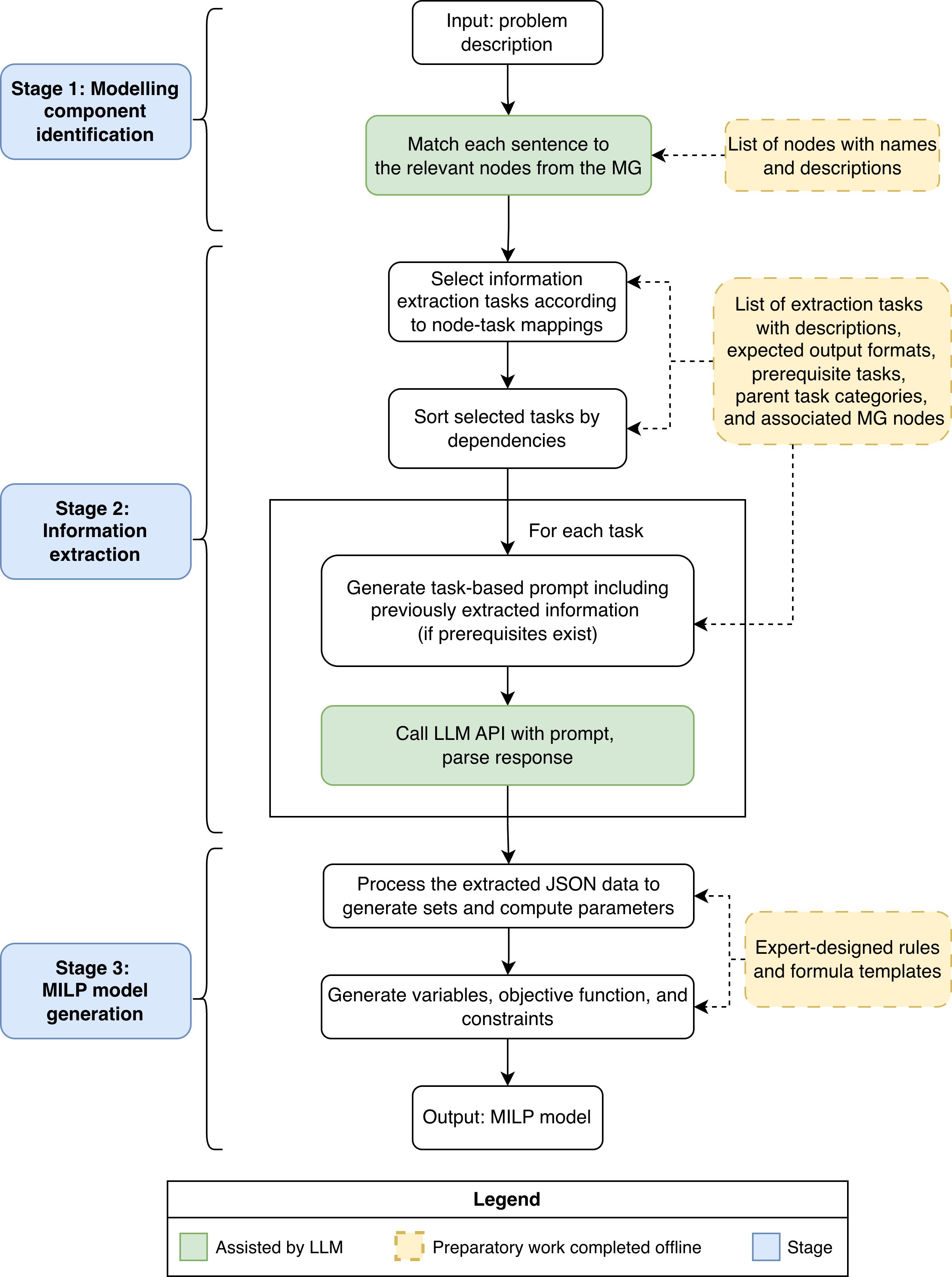}
  \caption{Main steps of the automatic formulation framework}
  \label{fig:flow1}
\end{figure}

By integrating semantic matching, structured knowledge representation, and controlled model synthesis, the framework provides a reproducible, interpretable, modular approach to automated MILP formulation. The main workflow of the proposed framework is illustrated in \autoref{fig:flow1}, while a step-by-step example, including input, transformation stages, and output, is provided in \autoref{fig:flow2}.

\subsection{Modelling graphs}

To ensure accurate and structured MILP model generation, our framework incorporates MGs as domain-specific knowledge modules. These graphs are designed to capture the fundamental components, including entities, parameters, decision variables, objectives, and constraints, required for formulating specific optimisation problem types, e.g., workforce scheduling problems. For each considered problem type, we develop a dedicated MG; examples include shift scheduling (see \autoref{fig:shiftKG} and \autoref{fig:shiftKG2}) and days-off scheduling (see \autoref{fig:daysoffKG} and \autoref{fig:daysoffKG2}). 
Each MG is further transformed into three distinct representations that enable automated retrieval across the stages of the framework: node descriptions, information extraction tasks, and mathematical formula templates. These three resources are pre-constructed during the framework development phase. Their specific roles will be elaborated in the following sections.

\begin{landscape}
\begin{figure}[!ht]
  \centering
  \includegraphics[width=0.8\linewidth]{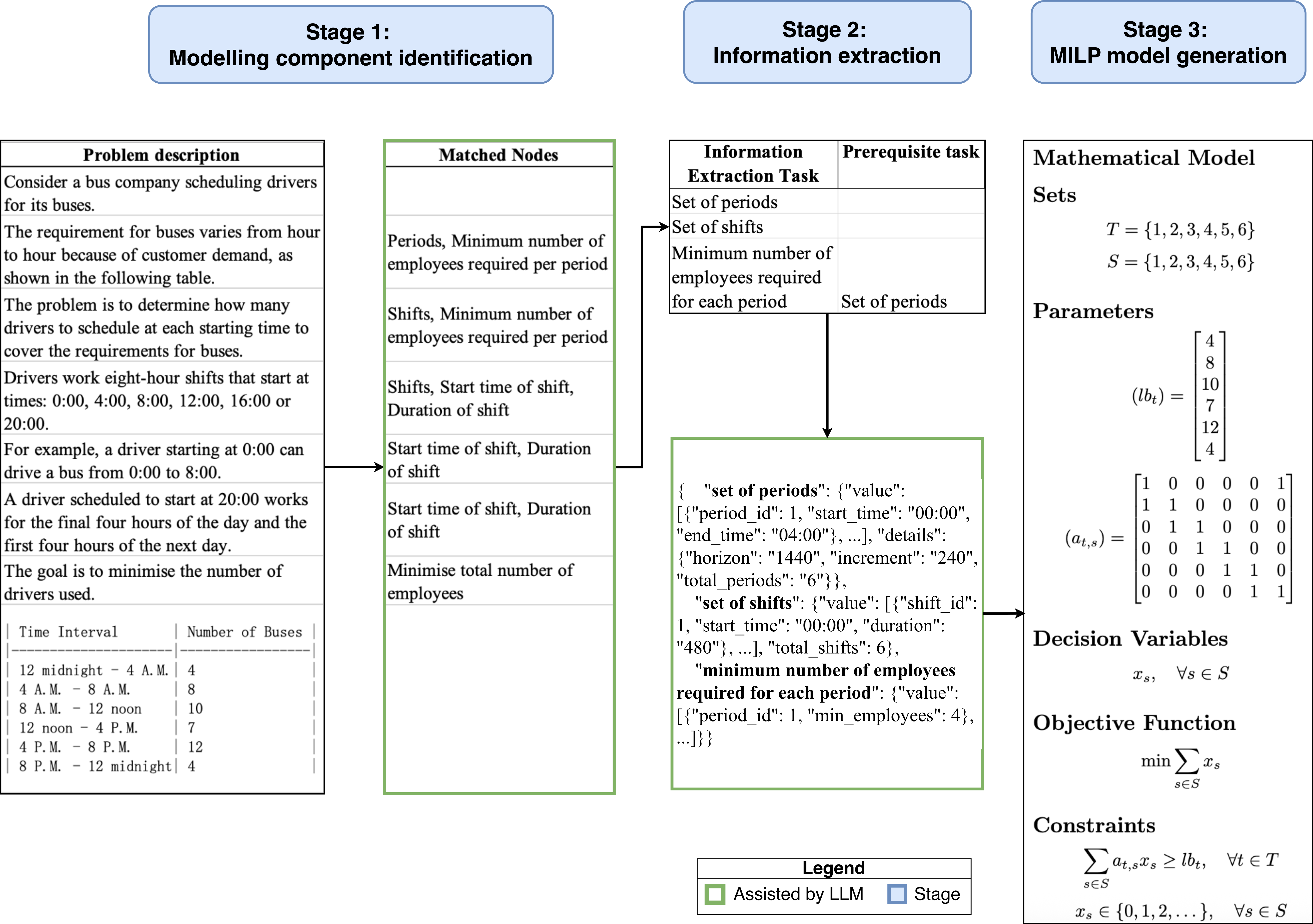}
  \caption{Implementing the automatic formulation framework on a problem instance of shift scheduling from our test set}
  \label{fig:flow2}
\end{figure}
\end{landscape}

\begin{figure}[!ht]
  \centering
  \includegraphics[width=1\linewidth]{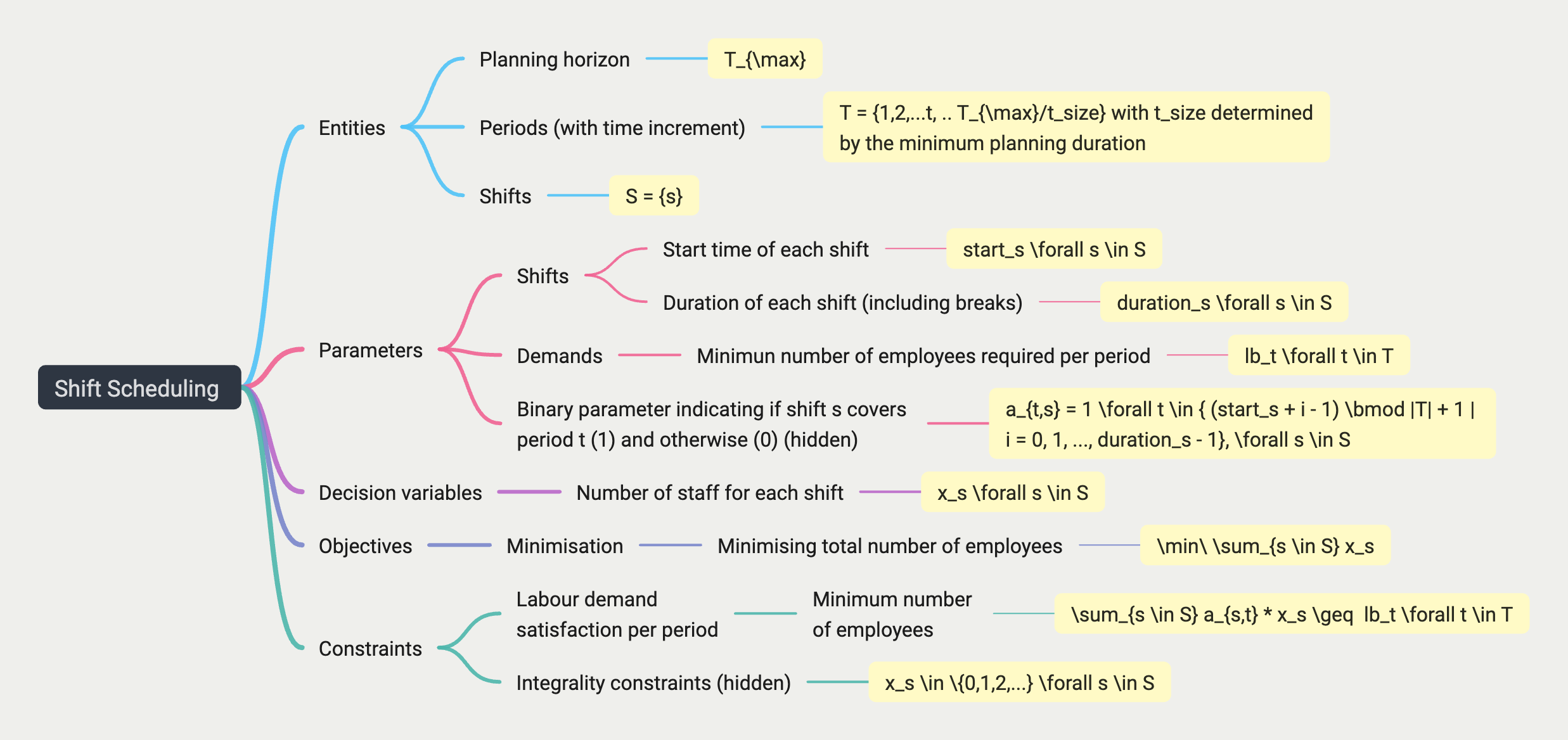}
  \caption{Part of the shift scheduling MG}
  \label{fig:shiftKGsnapshot}
\end{figure}

We use a simplified version of the shift scheduling MG (see \autoref{fig:shiftKGsnapshot}) to introduce the core ideas and the fundamental components. \autoref{fig:shiftKGsnapshot} only includes the essential nodes for the simplest instances of this problem type, such as the instance shown in \autoref{fig:flow2}. The full MG for shift scheduling, presented in \autoref{fig:shiftKG}, incorporates additional elements, including the scheduling of breaks and overtime, as well as more types of constraints and objective functions.

Each MG consists of several interconnected layers that define the components of an optimisation problem:
\begin{itemize}
    \item \textbf{Entities}: The MG encodes fundamental scheduling elements such as time periods, the planning horizon, and shifts. The planning horizon is denoted by $T_{\max}$, while time is discretised into a set of periods $T = \{1, 2, \dots, T_{\max}/t_{\text{size}}\}$, based on the minimum planning duration $t_{\text {size}}$. The set of shifts denoted by $S = \{s\}$ serves as the basis for workforce allocation. Each shift is specified by its start time and duration, which are provided as parameters.
    
    \item \textbf{Parameters}: The MG defines essential numerical values required for the model formulation. These include shift attributes (i.e., the start time of each shift, $\text{start}_s, \forall s \in S$, and the duration of each shift, $\text{duration}_s, \forall s \in S$), and lower bounds involved in labour demand constraints (e.g., minimum number of employees per period, $\text{lb}_t, \forall t \in T$). The binary parameter $a_{t,s}$, which indicates whether shift $s$ covers period $t$, is computed using predefined logic and is incorporated into the constraint formulations. This type of parameter is not explicitly stated in the problem description but must be derived based on the attributes of each shift. 

    \item \textbf{Decision variables}: The number of staff assigned to each shift is denoted by $x_s, \forall s \in S $, forming the basis for workforce scheduling decisions.  

    \item \textbf{Objectives}: The objective is usually to minimise the total number of employees required or the associated scheduling cost.  

    \item \textbf{Constraints}: The essential constraints include labour demand satisfaction constraints ensuring that staffing levels meet or exceed minimum requirements ($\sum_s a_{t,s} x_s \geq \text{lb}_t, \forall t \in T$) and non-negativity and integrality constraints enforcing non-negative and integer values for decision variables ($x_s \in \{0,1,2,\dots\}, \forall s \in S$). Identifying which active shifts cover the labour demand for a given period is both crucial and challenging, especially when the scheduling of in-shift breaks and overtime is involved. This process is central to modelling shift scheduling problems, as it requires additional reasoning and computation. Problem descriptions typically do not explicitly or straightforwardly specify such details. If multiple breaks and multiple break windows are involved, the formulation will be more complex.
\end{itemize}

By explicitly encoding the relationships between scheduling components, the MG serves as a guidance mechanism that compensates for the lack of optimisation expertise in LLMs. The MG ensures completeness by activating all relevant components based on both the textual problem description and predefined dependency and relevance rules. Instead of relying on the LLMs to infer the correct model structure, the MG systematically determines which components should be included, reducing the risk of missing constraints, parameters, or variables.  


\subsection{Stage 1: Modelling component identification}

Given a problem description for a specific problem type, the first stage of the modelling process involves identifying which components from the corresponding MG are relevant. To this end, an LLM is employed to segment the text into individual sentences and match each sentence to up to three of the most relevant nodes in the MG.

As noted earlier, three knowledge resources developed offline are used to support the overall framework. The first of these is essential in this stage. 
Specifically, each MG node is supplemented with an extended description in addition to its name, clarifying its meaning within the context of mathematical optimisation.
These node names and descriptions are compiled into a reference list and incorporated into the LLM prompt, as shown in \autoref{lst:component}, enabling accurate alignment between natural language text and formal modelling components for the specific problem type.

\subsection{Stage 2: Information extraction via task-based prompting}

Once a set of relevant nodes has been identified, the framework proceeds to extract instance-specific details associated with the knowledge represented by those nodes. 
For example, if a node corresponds to a parameter, its concrete numerical value must be retrieved from the problem description. 
To achieve this, a set of information extraction tasks is defined for each modelling graph, covering its set and parameter nodes. As each modelling graph is constructed for a specific optimisation problem type, the corresponding task sets are tailored accordingly.

Notice that our framework does not rely on LLMs to identify decision variables or to formulate objective functions and constraints. For complex optimisation problems, these modelling tasks require specialised techniques and domain knowledge, which LLMs sometimes struggle to apply correctly and consistently. Instead, such responsibilities are handled by the framework's rule-based logic, which is driven by expert knowledge encoded in the modelling graph.

After selecting and ordering the applicable tasks, the framework executes them sequentially using an LLM. 
The second preparatory resource supports this stage by providing a repository of prompt templates for information extraction tasks, as well as definitions of both the relationships between tasks and the mappings between tasks and their corresponding modelling graph nodes.
Each task entry specifies the task name, task description, expected output format, prerequisite tasks, parent task category, and the associated modelling graph node(s).

The relationships between tasks fall into two categories: prerequisite and hierarchical. Prerequisite relationships determine the dependency structure among tasks and are used to pass known information and play a central role in task selection and ordering. Hierarchical relationships are used to eliminate redundant or overlapping tasks during task selection by enforcing preference for more specific or fine-grained tasks. The full details of task selection and ordering will be discussed in later subsections.

For each task, the prompt given to the LLM includes the task name, task description, expected output format, the full problem description, and known information derived from its prerequisite tasks. This known information is provided in JSON format. The LLM output is parsed and stored as JSON, ensuring only valid structured data is retained.
\autoref{lst:min_number} shows an example of a prompt designed to extract the minimum number of employees required for each period. As this task depends on the definition of the set of periods, the prompt incorporates information extracted from the prerequisite task concerning period definitions as known input.
Additional prompt examples are provided in \ref{app:prompt}.

\subsubsection{Task selection via node activation}

Based on the nodes identified in Stage 1, the framework must determine which information extraction tasks should be executed. To enable this, a predefined mapping is established between MG nodes and extraction tasks, where each task is associated with one or more nodes. A task is selected if at least one of its associated nodes is activated. In this way, activating a node can trigger the activation of related nodes, producing a propagation effect across the modelling graph.
This mechanism serves two purposes. First, it mitigates the impact of omissions in the component identification stage by introducing flexibility and tolerance to partial recognition. 
Second, rather than designing a separate task for each node, the framework groups multiple related nodes into a single task and extracts related information together. This reduces redundancy, avoids duplicated effort, and improves the LLM's ability to semantically interpret the problem context and the task requirement, thereby enhancing the accuracy of information extraction.

If a task is selected, all of its prerequisite tasks are also selected. This further contributes to the robustness of the framework against incomplete or partially correct initial node identification results. Task selection thus propagates along the dependency structure, ensuring that all required preceding tasks are included in the extraction sequence.
In addition, if a selected task has an associated generalised task, which refers to a coarser-grained alternative covering the same modelling component, and both tasks appear in the execution list, then the generalised task is removed. This ensures that more specific and fine-grained tasks are prioritised, as they are linked to predefined modelling rules embedded in the framework. 
For example, instead of extracting the shift costs directly, the framework may extract per-period costs, which differ depending on whether the work is regular or overtime. The shift costs will then be calculated in our system based on the duration of the shifts and the costs of the periods. This approach improves both accuracy and consistency, as the parameter calculation is handled through rule-based logic rather than relying on the LLM's interpretation.

By carrying out parameter calculation and formula construction using the internal logic of the framework instead of through the LLM, the burden on the LLM is significantly reduced. Expert knowledge, represented through parameter derivation rules and formula templates for different types of constraints and objective functions, plays a central role in guiding the formulation process. This approach enhances the correctness and reproducibility of the generated models while offering greater control and interpretability in the formula construction process.

\subsubsection{Task dependency sorting via topological ordering}

Information extraction tasks are interdependent, as certain tasks require information extracted by others and must therefore be executed in a specific sequence. To resolve these dependencies, SMILO employs a topological sorting algorithm to determine a valid execution order.
A task dependency graph is first constructed, where each task is represented as a node, and a directed edge from task $T_i$ to task $T_j$ indicates that $T_j$ depends on the output of $T_i$, that is, $T_i$ is a prerequisite task for $T_j$. This graph structure enables the identification of an appropriate task sequence through topological sorting.

The sorting procedure begins by computing the in-degree of each node, which corresponds to the number of other tasks that must be completed before it can be executed. Tasks with zero in-degree, meaning they have no unmet dependencies, are placed into an initial queue. The framework then iteratively processes this queue by removing tasks, reducing the in-degree of their dependent tasks, and appending any new zero in-degree tasks to the queue. This process continues until all tasks have been ordered. The result is a task sequence in which each task is positioned only after all its prerequisites have been satisfied.

\subsection{Stage 3: MILP model generation} 

Based on the information extracted in Stage 2, the framework proceeds to generate the sets and compute the parameters required for MILP model construction, following predefined rules. 
Based on the generated sets, the framework defines the relevant decision variables.
For example, given a shift scheduling instance, Stage 2 extracts the necessary information to construct two fundamental sets: the set of periods and the set of shifts. Period information includes the start and end times of each period in 24-hour format, with each period assigned a unique numeric identifier $t$, forming the set $T$. Shift information specifies the start time and duration of each possible shift, each assigned an identifier $s$, forming the set $S$.
Based on the constructed set $S$, the framework defines the decision variables $x_s$, $\forall s \in S$, representing the number of staff assigned to shift $s$.
A key computational component in this stage is the binary parameter matrix $a_{t,s}$, which encodes whether shift $s$ covers period $t$. This matrix is derived from the extracted start time and duration of each shift, with appropriate handling of cyclic scheduling, where shifts may extend across the boundary of the planning horizon. If overtime is permitted, the framework additionally computes the binary parameter matrix $v_{t,o,s}$, which specifies whether shift $s$ extended by $o$ overtime periods covers time period $t$. Similarly, if breaks are allowed, other binary parameters are computed in a similar way.

Once all extracted information has been processed and all derived parameters have been computed, the framework automatically assembles a complete MILP model in LaTeX format. This model is generated by combining the predefined mathematical notation templates associated with each modelling component. The resulting formulation is structured in a standard format, beginning with the sets, parameters, and decision variables, followed by the objective function and the constraints, as illustrated in \autoref{fig:flow2}.

The third preparatory resource supports this stage by defining the rules for parameter computation, the logic for variable definition, and the formula templates for various types of constraints and objective functions. For example, the labour demand constraint in the complete modelling graph (see \autoref{fig:shiftKG}) is specified in its most general form, accounting for both breaks and overtime. When the framework detects that breaks and/or overtime are not present in a particular instance, it provides a simplified formulation of the constraint, omitting the corresponding parameters and variables.

\section{Experimental Evaluation} \label{sec:experimental} 

In this section, we present a systematic evaluation of our automated MILP model generation framework. 
The experiment aims to assess the framework's ability to correctly process natural language problem descriptions, generate accurate MILP models and solve them successfully. 
The evaluation is based on 20 test instances, each corresponding to a distinct natural language problem description. We first describe the test instances and their characteristics, followed by the evaluation metrics used to assess the framework's performance. Then we present and discuss the experimental results.  

\subsection{Setup and test instances}  

We evaluated both our proposed multi-stage framework (with GPT-4o-2024-08-06 as the deployed LLM) and a baseline approach GPT-4o-2024-08-06. For each of the 20 test instances, we conducted five independent trials with each method and measured the correctness and consistency in generating mathematical formulations. 
Hereafter, GPT-4o refers specifically to version GPT-4o-2024-08-06. In all experiments, the temperature parameter of GPT-4o was set to zero to ensure deterministic outputs, while all other parameters were kept at their default values.

The proposed framework incorporates GPT-4o in two stages: component identification and information extraction. Based on the extracted information, the final mathematical model is produced by a predefined code implementation rather than by the large language model itself.

To enable automatic solution and evaluation, we developed a code generation prompt (see \autoref{lst:code1} and \autoref{lst:code2}) that guides GPT-4o to convert the LaTeX-formatted mathematical model into executable Python code, which can be solved using the Gurobi solver. This conversion procedure was applied consistently across all experimental settings.

As a baseline, we used GPT-4o in a zero-shot setting. In this setting, the large language model is directly prompted to generate LaTeX-formatted mathematical models from natural language problem descriptions. The resulting models are subsequently converted into Python code using the same code generation prompt.
We designed two modelling prompts: one with minimal guidance, referred to as \textit{the simple prompt} (see \autoref{lst:gpt1}), and another with more detailed and rigorous instructions, referred to as \textit{the standard prompt} (see \autoref{lst:gpt2}). \smallskip

We developed a novel dataset of 20 problem instances, consisting of 10 shift scheduling problems and 10 days-off scheduling problems. 
These two groups of instances were designed to evaluate two variants of our framework, which are built on the same architecture but adapted to different problem types via problem-type-specific modelling graphs.
Specifically, the shift scheduling instances were used to test the framework integrated with the shift scheduling modelling graph, while the days-off scheduling instances were used to evaluate the framework incorporating the days-off scheduling modelling graph.

The dataset was constructed to reflect a diverse range of modelling scenarios by combining various attributes commonly encountered in these two types of decision-making problems. Each instance includes a natural language problem description, often with tables, a reference mathematical model in LaTeX, and a known optimal objective value (obtained by solving the Python code with Gurobi). All instances are feasible, and their corresponding models admit optimal solutions. 

Verifying the optimal value enables a fully automated evaluation approach. 
Compared with automatic evaluation based on direct comparison against ground truth formulations, this method offers greater flexibility by allowing multiple valid formulations of the same problem. 
A given problem instance in our dataset may admit more than one correct formulation. For reference, we provided one mathematical model for each instance to facilitate the use of the dataset by other researchers. However, using a different formulation does not necessarily indicate an error. 
In the automatic evaluation, we compared the optimal value obtained by solving the generated model with the known optimal value from the dataset. 

To enhance the reliability of the evaluation, we also conducted a round of manual verification of the mathematical models. Specifically, the 20 models generated by our framework and the 20 models generated by the baseline approach using the standard prompt in the first trial were reviewed by human experts. This verification process assessed whether each model was fully correct when the optimal value matched the reference, and also identified reasons for incorrect solutions. Importantly, valid models with alternative formulation from the reference would not be marked as incorrect by the human experts.

Readers may find examples of the prompts given to the baseline models as well as those used in the multi-stage framework in \ref{app:prompt}. 
Data used in this paper is available at: \url{https://anonymous.4open.science/r/C28M9232F635/}.

\subsection{Evaluation metrics} 

The performance of each method was assessed using two metrics:  
\begin{itemize}
    \item Model Accuracy ($\rm MA$): Measures whether the generated MILP model is \textit{fully correct}. A model is considered correct only if all components (sets, variables, parameters, constraints, and the objective function) are correctly defined and represented.
    \item Execution Accuracy ($\rm EA$): Measures the correctness of the optimal value obtained by solving the generated model. An instance is considered correctly solved if the computed optimal value exactly matches the provided ground truth. 
\end{itemize}

Each metric is calculated as the proportion of correctly processed instances out of the 20 total test instances. The metric $\rm MA$ reports the results of manual verification of the mathematical models generated in the first trial. The metric $\rm EA$ is computed over all five trials, along with its average. The formal mathematical formulas of the two metrics are given below:

\[
\rm MA = \frac{\text{Number of test instances with fully correct MILP models}}{\text{Total number of test instances}},
\]

\[
\rm EA = \frac{\text{Number of test instances with correct optimal values}}{\text{Total number of test instances}}.
\]

\subsection{Experimental results}

\subsubsection{Overall evaluation of framework and baseline} 

We evaluated the proposed multi-stage framework SMILO on the test set and compared its performance with the baseline method. Both GPT-4o (zero-shot setting) and the multi-stage framework with GPT-4o were tested on all instances, with each instance repeated five times, resulting in 100 generated models per method. The correctness of the MILP models produced by each method was evaluated automatically, based on whether the generated model yielded the correct optimal value. In addition, the models generated in the first trial were manually reviewed by human experts. The findings from manual inspection are discussed in the next subsection.

As shown in \autoref{tab:ea_ma_comparison}, the multi-stage framework, which integrates knowledge representations with the large language model, significantly outperforms the one-step model generation approach in terms of both automatic and manual evaluation metrics. The multi-stage framework achieves an average Execution Accuracy (EA) of 90\% across five trials and a model accuracy (MA) of 90\% in the first trial, demonstrating strong reliability and reproducibility in producing correct and consistent models. For two problem instances, the framework consistently produced incorrect models across all five trials. These failures were caused by the same type of modelling error, which we analyse in the next subsection.

\begin{table}[htbp]
\centering
\caption{Comparison of execution accuracy (EA) across five trials and model accuracy (MA) in Trial 1}
\label{tab:ea_ma_comparison}
\begin{tabular}{|p{5.15cm}|*{6}{p{1.2cm}}|p{1.2cm}|}
\hline
\textbf{Method} & \multicolumn{6}{c|}{\textbf{EA}} & \textbf{MA} \\
\cline{2-8}
& Trial 1 & Trial 2 & Trial 3 & Trial 4 & Trial 5 & Avg & Trial 1 \\
\hline
SMILO & 90\% & 90\% & 90\% & 90\% & 90\% & 90\% & 90\% \\
GPT-4o with standard prompt & 35\% & 35\% & 35\% & 55\% & 50\% & 42\% & 35\% \\
GPT-4o with simple prompt & 35\% & 35\% & 35\% & 35\% & 55\% & 39\% & -- \\
\hline
\end{tabular}
\end{table}

Regarding the performance of the baseline approach, using the standard prompt led to a modest improvement in average EA, increasing from 39\% with the simple prompt to 42\%. However, none of the models generated using the simple prompt resulted in infeasible solutions, whereas the standard prompt led to 7 infeasible models out of 100. This indicates that while the standard prompt may guide the model toward higher overall accuracy, it also introduces a higher risk of generating infeasible models.

GPT-4o's behaviour exhibited substantial variability. For instance, under the standard prompt setting, it produced correct models for 9 problem instances in some trials but generated incorrect models for the same instances in others. For 8 problem instances, the model consistently failed across all five trials, although the types of errors varied. In contrast, three problem instances were consistently solved correctly in every trial. They were the simplest problems in our dataset.
The fluctuation in correctness across trials for the same instance suggests a lack of output stability, even under deterministic decoding settings. The consistently failed cases, despite diverse error types, indicate that certain problem attributes remain challenging for zero-shot prompting. Conversely, the consistently successful cases suggest that for specific types of simple problems, GPT-4o is capable of reliably applying the appropriate modelling logic.

\begin{figure}[!ht]
  \centering
  \includegraphics[width=1\linewidth]{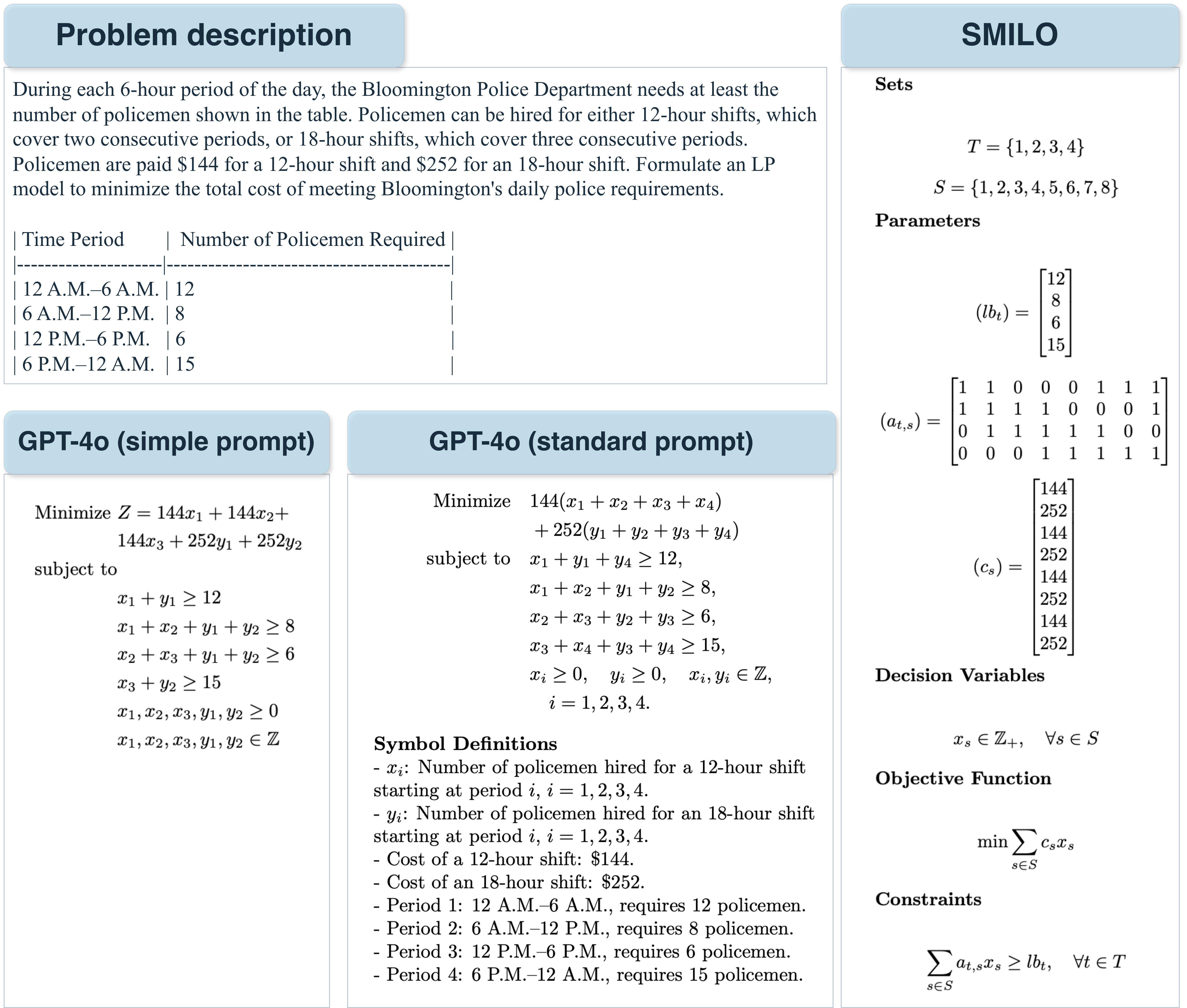}
  \caption{Comparison of MILP models generated by GPT-4o using the simple prompt, GPT-4o using the standard prompt, and the proposed SMILO framework. GPT-4o using the simple prompt misdefined the decision variables. GPT-4o using the standard prompt generated incorrect constraints. SMILO generated a correct model.}
  \label{fig:model_comparison}
\end{figure}

To illustrate the differences in model quality between the baseline under the two prompting strategies and our proposed framework, we present an empirical example comparing the generated mathematical formulations, as shown in \autoref{fig:model_comparison}. The example is a shift scheduling problem instance from our dataset, and the models were produced during the first evaluation trial. GPT-4o with the simple prompt incorrectly defined the decision variables by failing to identify all possible shifts. GPT-4o with the standard prompt correctly defined the variables but misidentified the coverage of shifts over periods, resulting in incorrect formulations of the first four constraints. SMILO generated a complete and correct model.

Next, we shall look into the cases when incorrect models were generated, providing a detailed analysis of the associated error types.

\subsubsection{Analysis of error cases}

To analyse the types of modelling errors, we manually reviewed the 20 models generated in the first trial by our proposed framework and the 20 models generated by the baseline approach using the standard prompt. 
The review showed that whenever the optimal value was correct, the corresponding mathematical model (in LaTeX format) was also correct. 
In cases where the optimal value was incorrect, the underlying issue was generally related to model formulation rather than to code generation or execution. Among the 40 reviewed models, only one exception was found in which the mathematical model was correct, but an error in the generated code altered the meaning of the model and led to an incorrect optimal value.
Therefore, our error analysis focuses exclusively on modelling-related errors. The results reveal that the primary challenges lie in correctly identifying the shift types, defining appropriate decision variables, and handling their coverage over time periods.
We classified the modelling errors produced by the baseline method and our framework into nine categories, as summarised in \autoref{tab:error_comparison}. 

\begin{table}[htbp]
    \centering
    \caption{Comparison of error types between the baseline and SMILO}
    \label{tab:error_comparison}
    \begin{tabular}{|p{5cm}|p{1.7cm}p{1.3cm}p{1.3cm}|p{1.7cm}p{1.1cm}p{1.1cm}|}
    \hline
    \textbf{Error type} & \multicolumn{3}{c|}{\textbf{GPT-4o with standard prompt}} & \multicolumn{3}{c|}{\textbf{SMILO}} \\
    \cline{2-7}
    & \textbf{Number of errors} & \textbf{Overall \%} & \textbf{Error \%} & \textbf{Number of errors} & \textbf{Overall \%} & \textbf{Error \%} \\
    \hline
    Incorrect information extraction & - & - & - & 2 & 10\% & 66.67\%  \\
    Missing equality constraint & - & - & - & 1 & 5\% & 33.33\%  \\
    Incorrect coefficients in labour demand constraints & 8 & 40\% & 44.44\% & - & - & - \\
    Misinterpreted objective or constraint semantics & 3 & 15\% & 16.67\% & - & - & - \\
    Incorrect handling of cyclic horizon & 2 & 10\% & 11.11\% & - & - & - \\
    Uninstantiated symbolic parameters & 2 & 10\% & 11.11\% & - & - & - \\
    Miscomputed objective coefficients & 1 & 5\% & 5.56\% & - & - & - \\
    Missing break within shift & 1 & 5\% & 5.56\% & - & - & - \\
    Incorrect linking constraint & 1 & 5\% & 5.56\% & - & - & - \\
    \hline
    \end{tabular}
\end{table}

\subsubsection*{Cases missed by the multi-stage framework}

Among the models generated by the multi-stage framework, two types of modelling errors were observed. In the first trial, the models generated for two of the 20 instances were incorrect, both belonging to the days-off scheduling problem set. Across all five trials, these two instances consistently resulted in incorrect models. Further manual inspection revealed that the modelling errors were consistently caused by inaccurate information extraction in Stage 2 of the framework. In both cases, the component identification in Stage 1 was consistently correct and successfully triggered the intended information extraction tasks. 
It is also worth noting that the baseline model failed to correctly formulate the first instance in all five trials, while for the second instance, it succeeded in four out of five trials.

In the first problematic instance, the objective was to minimise weekend bonuses. Employees who work on Saturdays and Sundays receive a bonus of \$25 per day. This implies that the cost of a shift depends on whether it includes work on a weekend day. However, in all five runs, the extracted shift costs were incorrect, leading to erroneous coefficients in the objective function. 
The second instance included a constraint requiring that at least half of the employees have weekends off, meaning they do not work on both Saturday and Sunday. The intended mathematical formulation ensures that the total number of employees who rest on Saturday and Sunday is greater than or equal to half of the total number of employees. However, in all five trials, the left-hand side of the inequality was incorrectly constructed using variables that represent the number of employees whose shifts start on Saturday and Sunday, rather than those who are actually off on those days.
It is important to note that, for both information extraction tasks, the prompts provided fully instantiated and accurate shift set information derived from the corresponding problem instances.
Each shift was specified by its identifier, start day, and duration.
With this information, it is straightforward to infer the working and resting days of each shift over the weekly planning horizon. Therefore, these errors fundamentally stem from GPT-4o's failure to reason about which shifts are active on specific days, such as Saturday or Sunday.

After completing the main experiments, we conducted additional explorations to examine whether improved prompting could mitigate such errors. Taking the first instance as an example, we modified the prompt by explicitly providing each shift's work days and rest days in natural language format (e.g., \texttt{"rest\_days": ["Saturday", "Sunday"]}), instead of providing its start day and duration. Despite this more intuitive representation, GPT-4o's responses remained incorrect. For example, a shift starting on Saturday and spanning five consecutive days was still incorrectly calculated to receive a weekend bonus of \$25, although it should have received \$50 in total.

In the first instance, there was also one equality constraint that was correctly identified in Stage 1 of all five trials, but for which the extraction of Stage 2 returned no information. In other words, GPT-4o failed to recognise this constraint in the problem description and treated it as nonexistent. The intended constraint specifies that the total number of employees assigned to all shift types should equal 25, ensuring that no employee remains unassigned. This constraint should be formulated as an equality, where the sum of all shift-assignment variables equals a fixed constant.
This failure may be attributed to GPT-4o's insensitivity to less explicitly stated constraints in natural language. In this case, the original sentence in the problem description did not present the constraint in a direct or formulaic manner, which may have made it more difficult for the LLM to recognise it as an equality constraint.

Outside the main experiments, we made targeted efforts to improve the constraint recognition by enhancing the prompt. Specifically, we incorporated two illustrative examples, each consisting of a natural language description paired with the corresponding extracted constraint information, such as identifying a restricted subset of shifts and specifying an upper bound on the total number of employees assigned to them. These examples were included in the dedicated prompt designed for extracting the information required to formulate such equality constraints. However, even with this additional contextual guidance, GPT-4o consistently returned an empty output for this constraint type.

\subsubsection*{Cases missed by the baseline approach}

In the first trial, the baseline approach generated incorrect formulations for 13 out of the 20 problem instances. These errors can be grouped into seven categories based on their nature and source. The most frequent category involved incorrect coefficients in the labour demand constraints. In eight cases, these coefficients were incorrect due to misidentification of key shift attributes, including duration, start time, lunch break duration, and lunch break start time. These misinterpretations resulted in inaccurate inference of workforce availability across time periods, which is fundamental for constructing correct period-related labour demand constraints.

A second type of error, observed in three models, involved misinterpretation of the intended semantics of the objective function or individual constraints. These issues often stemmed from GPT-4o's inability to correctly infer which shifts included work on weekends and which did not, leading to errors in the formulation of constraints related to weekend bonuses or weekend rest requirements.

In two models, the cyclic nature of the planning horizon was not properly handled. Specifically, the formulations failed to recognise that shifts extending beyond the final period continue into the initial period of the subsequent cycle. This omission caused an underestimation of labour availability during the early periods of the schedule.

Two models failed to fully instantiate parameters as required by the prompt, that is, they presented symbolic representations without assigning specific numerical values. For instance, in one case, the LaTeX formulation included a binary coverage parameter $a_{s,k}$, representing whether shift $s$ covers period $k$, but omitted its concrete values from the model and the associated list of symbol definitions in the output. This omission shifted the responsibility for computing the binary matrix to the code generation stage. However, necessary contextual details, such as the start and end time of each period, were no longer present in the output, making accurate computation impossible, even though the shift start times and durations in hours were provided. 
This type of failure reflects a broader issue: symbolic notations are introduced without being accompanied by the necessary instance-specific data. As a result, the generated model lacks the information required to support downstream processes, such as code execution, undermining automation and correctness.

In addition to these primary error types, several less frequent issues were also observed. One model contained incorrect coefficients in the objective function due to miscalculated shift costs. Another model omitted the requirement of scheduling lunch breaks within shifts. A third model introduced errors in linking constraints between variables. 
A summary of the error types and their frequencies is presented in \autoref{tab:error_comparison}. In the table, the column ``Overall Percentage" shows the frequency of each error type relative to the total number of test instances (20 in total), while ``Error Percentage" indicates the proportion of each error type among all modelling errors identified in the first trial. Both metrics are calculated separately for each method.

\section{Conclusion} \label{sec:conclusion}

In this paper, we propose SMILO, an expert-knowledge-driven framework for generating MILP models from natural language problem descriptions.
SMILO adopts an interpretable three-stage pipeline, supported by problem-type-specific modelling graphs and their associated resources, including modelling component descriptions, task-specific prompt templates for information extraction, and expert-defined templates for variable definitions, parameter computations, and expression generation. 
To support evaluation, we construct a comprehensive dataset of shift scheduling and days-off scheduling problems, which reflect rich structural complexity such as diverse work patterns, break allocation, and overtime rules.
SMILO achieves 90\% accuracy across five evaluation trials, as evaluated by expert validation of model correctness and alignment with known optimal objective values, outperforming one-step LLM baselines by at least 35\%. 
These results demonstrate the effectiveness of combining expert knowledge with LLM capabilities for robust and reliable optimisation modelling.

In addition to overall accuracy, we conduct a detailed analysis of error types across generated models. 
Compared with the one-step LLM baseline using the standard prompt, SMILO exhibits significantly fewer modelling errors, with only two distinct error types observed. 
This highlights the ability of SMILO in enhancing interpretability and diagnosability by structurally decomposing the modelling process: LLMs are responsible only for information extraction, guided by knowledge-guided prompts, while the mathematical formulation is handled through deterministic rule-based templates. 
This clear separation improves interpretation, error traceability, and consistency in model generation.

In practical deployments, SMILO can be applied in production and logistics contexts such as factory line staffing, last-mile driver shifts, and hospital or inspection team scheduling. Practitioner requirements written in natural language are converted by the system into solver-ready MILP models and executable code, reducing modelling time and enabling rapid model updates when staffing policies change. The modelling graph can be maintained as a reusable asset that captures specific rules on shifts, breaks and overtime. A light human-in-the-loop check on model accuracy and execution accuracy provides a clear acceptance criterion prior to deployment. This practice improves planning speed and knowledge transfer, while alleviating additional analytical workload for expert teams.

Several promising directions can further extend this work. 
First, extending SMILO to other optimisation problem types, such as inventory control, vehicle routing and dispatching, and capacity-constrained production planning, requires the development of new modelling graphs. The modelling graphs could be learned automatically by curating a dataset of annotated MILP models from academic publications. 
Second, a hybrid formula generation mechanism could be explored, where expert-defined templates are used to programmatically formulate standard modelling components and simultaneously serve as prompts to guide LLMs in generating formulations for optimisation scenarios beyond the predefined modelling graph. This hybrid approach could increase the flexibility and generalisability of the framework.
Third, transforming the modelling graphs into knowledge graphs may allow a richer semantic representation of modelling components, including their attributes and interdependencies. This transformation could facilitate more effective retrieval and reuse of modelling knowledge, thereby enhancing the automation and adaptability of the model generation process.



\section*{Declaration of interests}

The authors declare that they have no conflicts of interest.

\bibliography{sample-base}

\newpage

\appendix

\section{Modelling graphs}

\begin{figure}[!ht]
  \centering
  \includegraphics[width=1\linewidth]{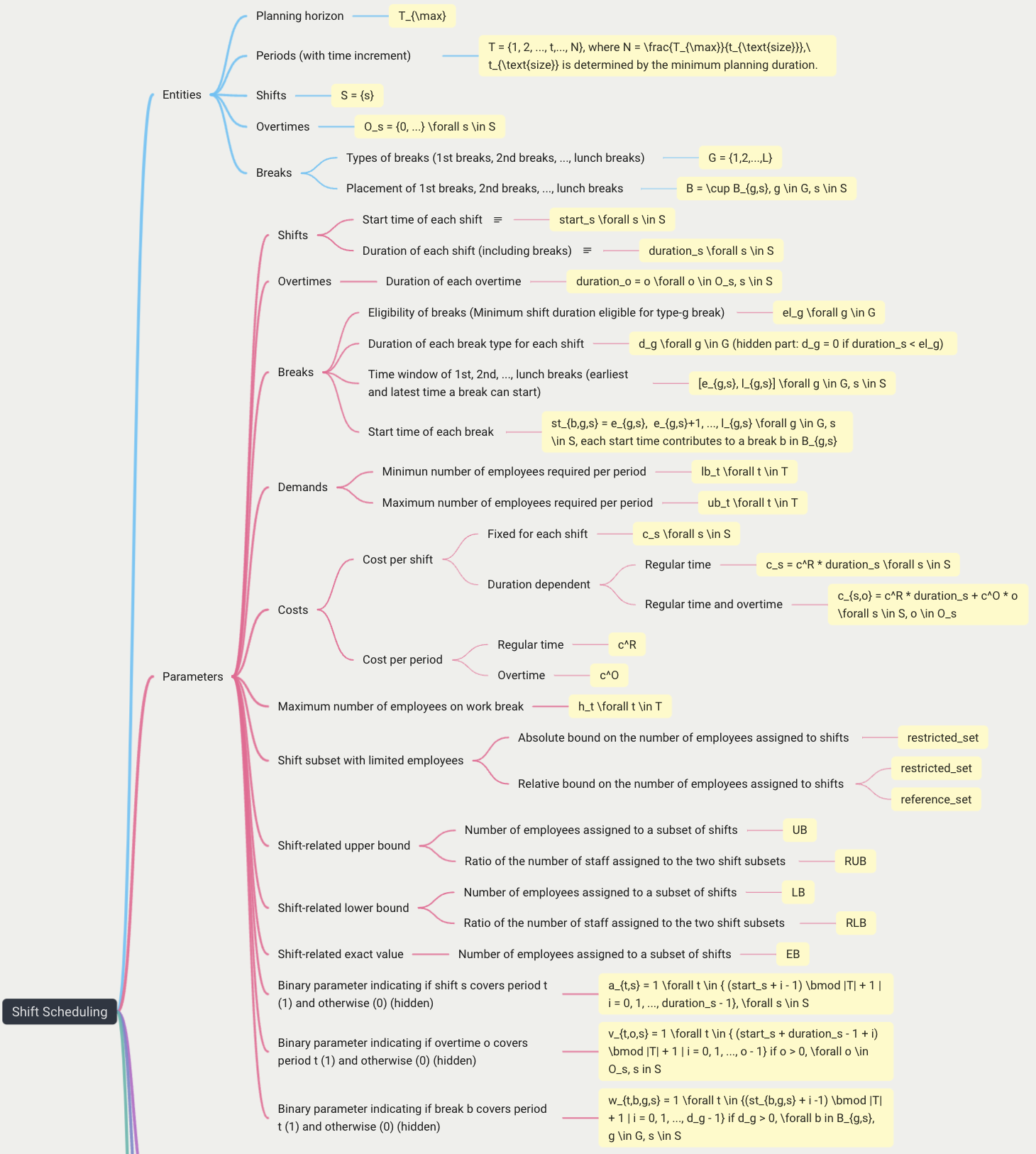}
  \caption{The modelling graph for shift scheduling (Part 1)}
  \label{fig:shiftKG}
\end{figure}

\begin{figure}[!ht]
  \centering
  \includegraphics[width=1\linewidth]{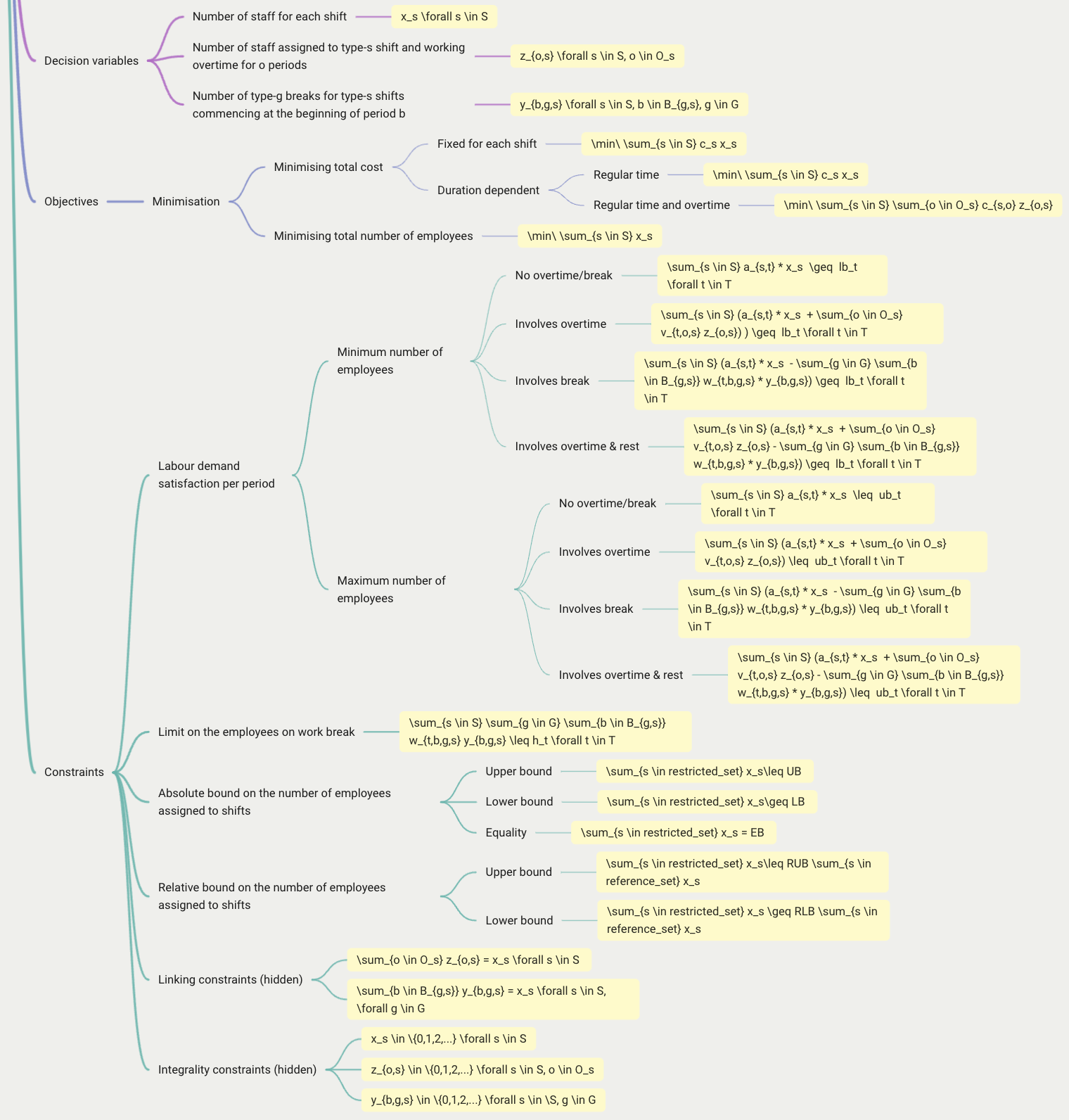}
  \caption{The modelling graph for shift scheduling (Part 2, continued)}
  \label{fig:shiftKG2}
\end{figure}

\begin{figure}[!ht]
  \centering
  \includegraphics[width=1\linewidth]{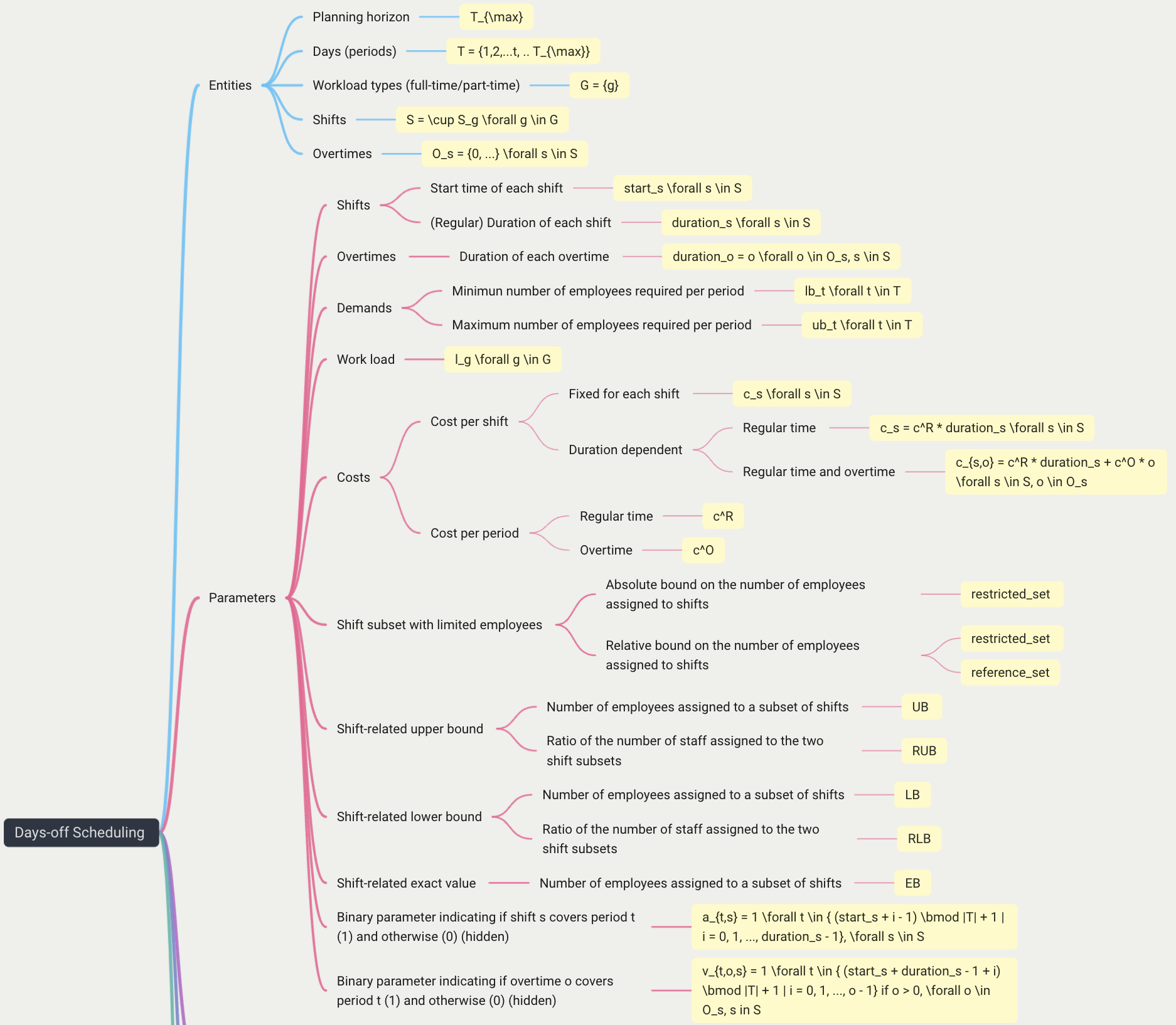}
  \caption{The modelling graph for days-off scheduling (Part 1)}
  \label{fig:daysoffKG}
\end{figure}

\begin{figure}[!ht]
  \centering
  \includegraphics[width=1\linewidth]{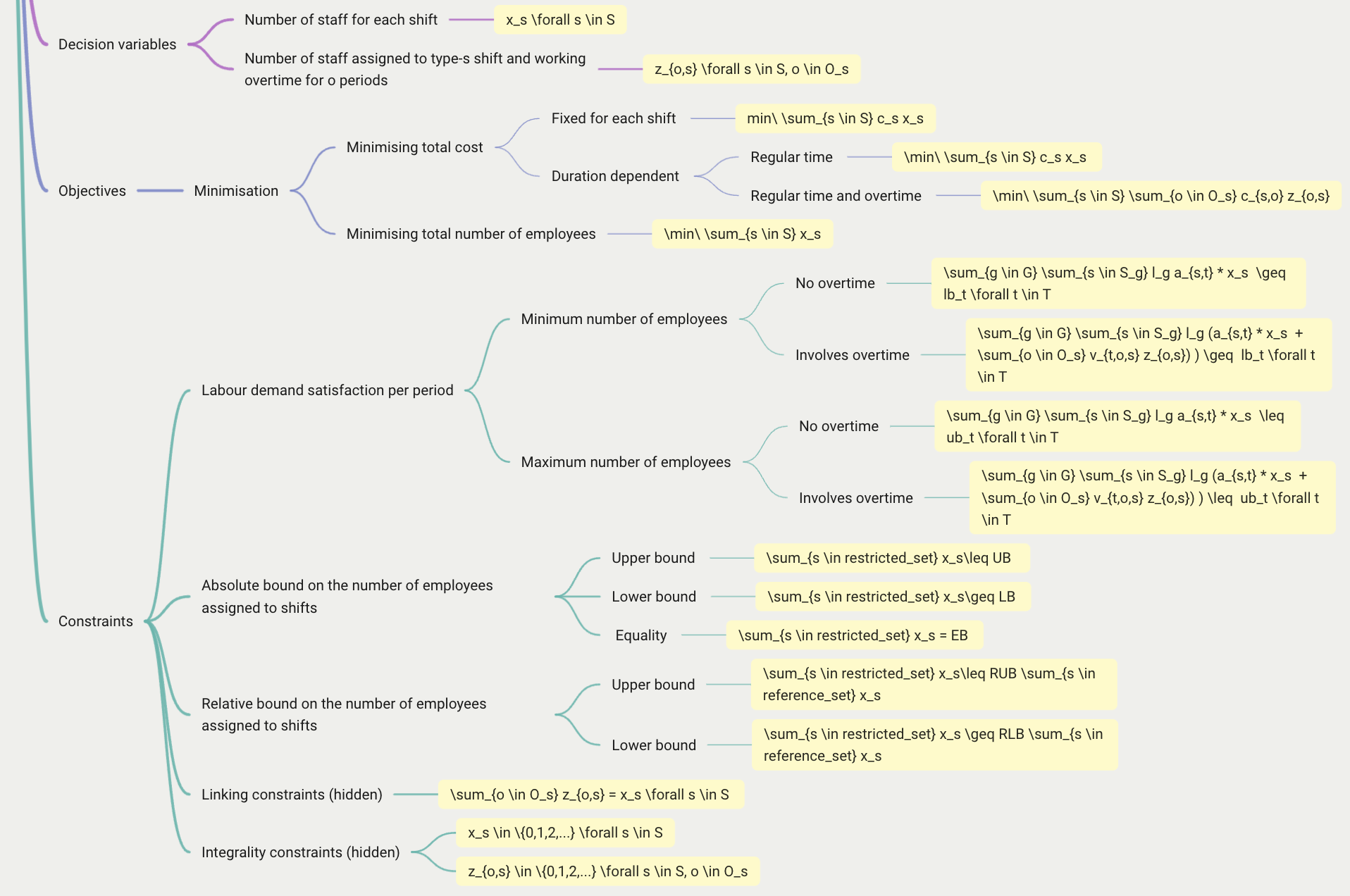}
  \caption{The modelling graph for days-off scheduling (Part 2, continued)}
  \label{fig:daysoffKG2}
\end{figure}

The shift scheduling MG is shown in \autoref{fig:shiftKG} and \autoref{fig:shiftKG2}. The days-off scheduling MG is shown in \autoref{fig:daysoffKG} and \autoref{fig:daysoffKG2}.

\section{Prompt template examples}
\label{app:prompt}
\subsection{Component identification prompt (SMILO)}

\lstset{
  basicstyle=\ttfamily\small,
  breaklines=true,
  frame=single,
  captionpos=b,
  columns=fullflexible,
  showspaces=false,
  showstringspaces=false,
}

\autoref{lst:component} shows the prompt template designed to identify relevant modelling components in the first stage of SMILO.
The placeholder \texttt{[Insert test instance here]} specifies the location at which the natural language problem description is dynamically inserted during runtime.
The node list shown here is a partial excerpt from the full node list of the shift scheduling MG used by SMILO. In practice, the full list includes all modelling components representing entities, parameters and objectives, along with their corresponding descriptions.
The ellipsis (\texttt{...}) marks the continuation of additional node entries not shown in this excerpt.

\begin{lstlisting}[caption={Component identification prompt used in SMILO}, label={lst:component}]
### Task: Identify Relevant Nodes for Each Sentence

#### Problem Description:
[Insert test instance here]

#### Modelling Graph Nodes:
Below are all the available modelling graph nodes. Each node has a **name** and a **description**. Your task is to match each sentence from the **Problem Description** to one or more of these nodes.

**Node Name**: Planning Horizon  
**Description**: The total time period over which scheduling or planning is conducted.

**Node Name**: Periods  
**Description**: The division of the planning horizon into smaller intervals, typically aligned with labour demand requirements. The increment determines the length of each period.

**Node Name**: Shifts  
**Description**: Represents all possible shift types involved in the problem. Each shift is defined by its start time and duration.

**Node Name**: Overtimes  
**Description**: Represents all possible overtime works involved in the problem.

**Node Name**: Breaks  
**Description**: The complete collection of all possible breaks across all types (e.g., lunch breaks, first breaks, second breaks). Each break is uniquely identified by its type, start time, and duration.

**Node Name**: Types of breaks  
**Description**: Different predefined types of breaks assigned to employees, such as first relief breaks, second relief breaks, and lunch breaks, based on specific rules or schedules.

...

#### Instructions:
1. **Extract all sentences from the Problem Description**, ensuring they are listed in the correct order.
2. **Assign a unique identifier to each sentence (e.g., "sentence_1", "sentence_2", etc.).**
3. **Ensure that each sentence is returned exactly as it appears in the original Problem Description.**
4. **Match each sentence to up to 3 most relevant nodes from the Modelling Graph.**
5. **Ignore tables as separate sentences.** If a sentence introduces a table, it must still be processed normally.
6. **Ensure the response follows the JSON format strictly** with the correct number of sentences.

#### Expected Output Format:
```json
{{
  "sentences": [
    "Sentence 1",
    "Sentence 2",
    "Sentence 3",
    ...
  ],
  "matches": {{
    "sentence_1": ["Node A", "Node B"],
    "sentence_2": ["Node C"],
    "sentence_3": []
  }}
}}

Return only the JSON response. Do not include explanations or additional formatting.
\end{lstlisting}

\subsection{Task-based prompts for information extraction (SMILO)}

\autoref{lst:min_number} shows an example of a task-based prompt designed to extract the minimum number of employees required for each period. 
In this example, the content under the headings ``Known Information" and ``Problem Description" is instance-specific and dynamically inserted at runtime.
This task depends on a prerequisite: extracting the set of periods. The prompt includes the extracted period information as known input.

\begin{lstlisting}[caption={Prompt for extracting the minimum number of employees required for each period}, label={lst:min_number}]
### Task Name
minimum number of employees required for each period

### Task Description
Extract the minimum number of employees required for each period in the set of periods. Follow these steps:

1. For each period, identify the associated demand in the problem description.
2. Output the results in the specified format, associating each period with its corresponding number of employees.

### Expected Output Format
{
  "value": [
    {
      "period_id": [integer],
      "min_employees": [integer]
    },
    ...
  ]
}

Please ensure your response strictly follows the **Expected Output Format** provided above. Avoid including any code implementations. The response must include a **valid JSON object** containing only the requested information.

### Known Information
{
    "set of periods": {
        "value": [
            {
                "period_id": 1,
                "start_time": "00:00",
                "end_time": "04:00"
            },
            {
                "period_id": 2,
                "start_time": "04:00",
                "end_time": "08:00"
            },
            {
                "period_id": 3,
                "start_time": "08:00",
                "end_time": "12:00"
            },
            {
                "period_id": 4,
                "start_time": "12:00",
                "end_time": "16:00"
            },
            {
                "period_id": 5,
                "start_time": "16:00",
                "end_time": "20:00"
            },
            {
                "period_id": 6,
                "start_time": "20:00",
                "end_time": "00:00"
            }
        ],
        "details": {
            "horizon": "1440",
            "increment": "240",
            "total_periods": "6"
        }
    }
}
### Problem Description
Consider a bus company scheduling drivers for its buses. The requirement for buses varies from hour to hour because of customer demand, as shown in the following table. 

The problem is to determine how many drivers to schedule at each starting time to cover the requirements for buses. Drivers work eight hour shifts that start at times: 0:00, 4:00, 8:00, 12:00, 16:00 or 20:00. For example, a driver starting at 0:00 can drive a bus from 0:00 to 8:00. A driver scheduled to start at 20:00 works for the final four hours of the day and the first four hours of the next day. The goal is to minimise the number of drivers used. 

| Time Interval       | Number of Buses |
|---------------------|---------------------------|
| 12 midnight - 4 A.M.| 4                        |
| 4 A.M. - 8 A.M.     | 8                         |
| 8 A.M. - 12 noon    | 10                         |
| 12 noon - 4 P.M.    | 7                         |
| 4 P.M. - 8 P.M.     | 12                         |
| 8 P.M. - 12 midnight| 4                         |
\end{lstlisting}

\autoref{lst:workload_shifts} presents the prompt designed to extract the set of workload-specific shifts in days-off scheduling problems. 
In this example, the content under the headings ``Known Information" and ``Problem Description" is instance-specific and dynamically inserted at runtime.
As this task depends on the definition of the set of days and the set of workload types, the prompt incorporates information extracted from the prerequisite tasks concerning period definitions and workload definitions as known input.
We also designed a more general prompt template for extracting the set of shifts, defined as the parent category task for extracting workload-specific shifts. Since extracting the set of workload-specific shifts is a more fine-grained task, it is executed in cases where both tasks are matched to a problem instance.

\begin{lstlisting}[caption={Prompt for extracting the set of workload-specific shifts}, label={lst:workload_shifts}]
### Task Name
set of workload-specific shifts

### Task Description
For each daily workload type, extract all valid shift patterns by identifying the possible start days and shift durations within the planning horizon. Follow these steps:

1. Extract the total number of consecutive workdays per shift for each workload type.
2. Determine possible start days:
   - If the problem description explicitly specifies valid start days, extract those.
   - Otherwise, assume that any day in the planning horizon is a valid start day.
3. Generate all valid shifts, where each shift is defined by:
   - Start day (must be a valid start day).
   - Duration (total number of consecutive workdays per shift).

### Expected Output Format
{
  "value": [
    {
      "workload_type": "[Workload type name]",
      "shifts": [
        {
          "shift_id": [integer],
          "start_day": "[Start day in planning horizon (e.g., Day 1, Day 2, ..., Day N)]",
          "duration": "[Number of consecutive working days (excluding overtime days)]"
        },
        ...
      ]
    },
    ...
  ],
  "total_shifts": [integer]
}

Please ensure your response strictly follows the **Expected Output Format** provided above. Avoid including any code implementations. The response must include a **valid JSON object** containing only the requested information.

### Known Information
{
    "set of periods": {
        "value": [
            {
                "period_id": 1,
                "day": "Day 1"
            },
            {
                "period_id": 2,
                "day": "Day 2"
            },
            {
                "period_id": 3,
                "day": "Day 3"
            },
            {
                "period_id": 4,
                "day": "Day 4"
            },
            {
                "period_id": 5,
                "day": "Day 5"
            },
            {
                "period_id": 6,
                "day": "Day 6"
            },
            {
                "period_id": 7,
                "day": "Day 7"
            }
        ],
        "details": {
            "horizon": "7 days",
            "increment": "1 day"
        }
    },
    "set of workload types": {
        "value": [
            {
                "workload_type": "Full-time",
                "workload": 1
            },
            {
                "workload_type": "Part-time",
                "workload": 0.5
            }
        ]
    }
}
### Problem Description
A post office requires different numbers of full-time employees and part-time employees on different days of the week. The table specifies the daily labour hours required for each day of the week. During each week, a full-time employee works 8 hours a day for five consecutive days, and a part-time employee works 4 hours a day for five consecutive days. These shifts can start on any day of the week.

The number of part-time employees should not exceed 25% of the total number of employees, including both full-time and part-time staff.

A full-time employee costs the post office $600 per week, whereas a part-time employee (with reduced fringe benefits) costs the post office only $200 per week. Formulate an LP model to minimise the weekly labour costs for the post office.

| Day of the week  | Number of hours required |
|-----------|--------------------------|
| Monday    | 136                      |
| Tuesday   | 104                      |
| Wednesday | 120                      |
| Thursday  | 152                      |
| Friday    | 112                      |
| Saturday  | 128                      |
| Sunday    | 88                       |

\end{lstlisting}

\subsection{LaTeX-to-Gurobi code generation prompt (all methods)}
To enable automatic solution and evaluation, we developed a code generation prompt that guides GPT-4o to convert the LaTeX-formatted mathematical model into executable Python code, which can be solved using the Gurobi solver. This conversion procedure was applied consistently across all experimental settings.
The placeholder \texttt{[Insert the Latex-formatted mathematical model of the test instance here]} specifies the location at which the Latex-formatted mathematical model is dynamically inserted during runtime.

\begin{lstlisting}[caption={System prompt for LaTeX-to-Gurobi conversion}, label={lst:code1}]
You are an expert in mathematical optimisation and Python programming using Gurobi.
Your job is to translate LaTeX-formatted mathematical programming models into Python code that can be executed with the Gurobi optimiser.            
Do not change any model structure or parameter values.
Use the `gurobipy` API. Ensure the script can run directly when copied into a .py file.
The final output must be a single Markdown code block starting with ```python and ending with ```.
Do not include any explanation or commentary before or after the code block.
You must retain Gurobi's default solver output from model.optimize() (do not disable it).
After solving the model:
- If an optimal solution is found, print the results in the following format:
- A line starting with: `Optimal objective value: ` followed by the numerical value
  - A line starting with: `Variable values:`
- Then print each variable that has a non-zero value in the format: `VarName = Value`
  - If no optimal solution is found, print: `No optimal solution found.`
  You may implement this logic however you see fit, but the printed output must follow this format exactly.
\end{lstlisting}

\begin{lstlisting}[caption={User prompt for LaTeX-to-Gurobi conversion}, label={lst:code2}]
Given the following mathematical programming model written in LaTeX, please convert it into executable Python code using the Gurobi optimiser.

Requirements:
- Do not change any of the model's sets, parameters, variables, or expressions. Follow the structure of the model closely. All variables should use the same names and indices as in the original model.
- Use the `gurobipy` API. The output must be a complete, standalone Python script that can run directly when saved as a .py file.
- Keep Gurobi's default solver output when calling model.optimize().
- After solving:
  - If an optimal solution is found, print the objective value using the format: `Optimal objective value: <value>`
  - Then print all variables with non-zero values, preceded by a line that says: `Variable values:`
  - If no optimal solution is found, print: `No optimal solution found.`
- Only output a single Markdown code block starting with ```python and ending with ```.
- Do not include any explanation before or after the code block.

Model:
[Insert the Latex-formatted mathematical model of the test instance here]
\end{lstlisting}

\subsection{Simple and standard prompts for one-step model generation (Baseline GPT-4o)}
As a baseline, we used GPT-4o to directly generate LaTeX-formatted mathematical models from natural language problem descriptions. 
We designed two modelling prompts: one with minimal guidance, referred to as \textit{the simple prompt}, and another with more detailed and rigorous instructions, referred to as \textit{the standard prompt}. \smallskip
The placeholder \texttt{[Insert the problem description of the test instance here]} specifies the location at which the natural language problem description is dynamically inserted during runtime.

\begin{lstlisting}[caption={The simple prompts for GPT-4o}, label={lst:gpt1}]
You are an expert in Mixed-Integer Linear Programming (MILP). When given a problem description, identify decision variables, objective functions, and constraints. Provide explanations first, and then output the complete mathematical model at the end in the following format:

### Begin MILP Model
... (LaTeX formatted model) ...
### End MILP Model

[Insert the problem description of the test instance here]

\end{lstlisting}

\begin{lstlisting}[caption={The standard prompts for GPT-4o}, label={lst:gpt2}]
You are an expert in Mixed-Integer Linear Programming (MILP). When given a problem description, generate the mathematical model by extracting decision variables, objective functions, and constraints.Follow these instructions carefully:

1. Decision Variables: Clearly define all decision variables, including their meaning, type (e.g., binary, integer, continuous), and associated indices. Indices must be derived from the problem description and explicitly defined with their domain.

2. Parameters and Constants: 
   - If the problem description contains specific numerical values for parameters (e.g., in tables or lists), you must not omit them. 
   - You may either:  
     (i) include the concrete values directly in the LaTeX model, or  
     (ii) use parameter symbols (e.g., \( b_i \)) but you must then include the full list of values (e.g., \( b_1 = 10, b_2 = 15, \ldots \)) in the Symbol Definitions section.
   - Do not use symbolic parameters without providing their corresponding values if those values are available in the input.

3. Intermediate Expressions: If a symbol is introduced as an expression derived from decision variables and/or known parameters, define it clearly as an intermediate expression.

4. Mathematical Model: Provide a complete and well-structured MILP model using LaTeX format. All symbols used in the model must be explicitly defined and, where applicable, linked to specific values. Avoid using undefined or abstract placeholders.

5. Output Structure:
   - First, explain the modelling process in plain language, including the meaning of the decision variables, the objective, and constraints.
   - Then output the mathematical model and symbol definitions in the following format:

### Begin MILP Model
... (LaTeX formatted model with objective, constraints, and all necessary terms. Use numerical values directly where appropriate) ...
### End MILP Model

### Begin Symbol Definitions
... (List of all symbols, including decision variables, parameters/constants with specific values if available, index sets, and intermediate expressions) ...
### End Symbol Definitions

Important: Never omit numerical values provided in the problem description. If a parameter symbol is used in the model, its values must be listed. Do not abstract concrete values unless necessary for clarity or scalability. 

[Insert the problem description of the test instance here]
\end{lstlisting}



\end{document}